\documentclass{article}[12pt]

\usepackage{mathrsfs}
\usepackage{amssymb}
\usepackage{amsmath}
\usepackage{amsbsy}
\usepackage{epsfig}
\usepackage{enumerate}
\usepackage{bm}

\newtheorem{thm}{Theorem}[section]

\newtheorem{cor}{Corollary}[section]
\newtheorem{rem}{\sc Remark}[section]

\def\proclaim#1{\par \smallskip\noindent {\bf #1}\bgroup\it\ }
\def\endproclaim{\egroup\par\smallskip}

 \newbox\TempBox \newbox\TempBoxA

\def\pr{\textsf{P} } % the symbol P for probability used the sans serif letter
\def\ep{\textsf{E} } % the symbol E for expectation used the sans serif letter
\def\Var{\textsf{Var}} % the symbol Var for covariance used the sans serif letter

\numberwithin{equation}{section}

\begin{document}

% "Title of the paper"
\title{\bf A Gaussian Process Approximation for a two-color Randomly Reinforced  Urns}
\date{}

\author{Li-Xin Zhang\footnote{Department of Mathematics, Zhejiang University, 310027, P.R. China}\\
 Zhejiang University} % Your postal address goes here.
\maketitle

\begin{abstract}
 We prove a Gaussian process approximation for the sequence of random compositions of a two-color randomly reinforced urn for both the cases with the equal and unequal reinforcement means. By using the Gaussian approximation, the law of the iterated logarithm and  the functional limit central limit theorem in both the stable convergence sense and the almost-sure conditional convergence sense are established. Also as a consequence, we are able to to prove that the distribution of the urn composition has no points masses both when the reinforcement means are equal and unequal  under the assumption of only finite $(2+\epsilon)$-th moments.


\end{abstract}

\section{Introduction}
\setcounter{equation}{0}
Asymptotic properties, including the strong consistency and asymptotic normality, of urn models and their applications   are widely studied  in recent years under various assumptions concerning the updating rules, for example, one may refer to Chauvin, Ponyanne and Sahnoun (2009), Bai, Hu and Rosenberger (2002), Bai and Hu (2005), Janson (2004,2006), Zhang, Chueng and Hu (2006) etc. In this paper, we consider a kind of two-color urn model, called the randomly reinforced urn (RRU) model, which is a generalization of  the original P\'olya urn (c.f., Eggenberger and P∩olya (1923), P\'olya (1931)). The main issue of this model different from most urn models in literature is that, as shown, the  proportions of balls in the urn will not converge to a non-extreme constant and the numbers of different type balls may increase in different speeds. This issue makes its asymptotic properties  quite different from the those of other urn models and difficult to study.

The RRU model is described as follows. Consider a two-color urn with  the initial  urn components  $\bm Y_0=(Y_{0,1}, Y_{0,2}))$, where $Y_{0,k}>0$ is the
number of type $k$ balls. The urn is sampled sequentially.   Suppose the urn components are $\bm
Y_m=(Y_{m,1},  Y_{m,2})$  after  $m$ samplings. At the $(m+1)$-th sampling,   a ball of type $k$ is drawn with  a probability
 $$p_{m+1,k}=\frac{Y_{m,k}}{|\bm Y_m|}, \;\; \text{ where }\; |\bm Y_m|=Y_{m,1}+ Y_{m,2}.$$
 And the sampled ball is replaced in the urn together with a nonnegative random number $U_{m+1,k}$ of balls of the same type $k$, generated from a distribution $\mu_k$ with mean $m_k>0$.   This  is the model introduced and formally named   the randomly reinforced urn in Mulier, Paganoni and Secchi (2006a). But it would appear   in earlier literatures in  different versions.    For example, Durham and Yu (1990) proposed a similar model  for sequential sampling in clinical trails.  In our RRU setting, the numbers of balls take positive real values, not necessary integers.
   When $U_{m+1,1}=U_{m+1,2}=\alpha$ is a constant and   a positive integer, a RRU is  the original P\'olya urn (c.f., Eggenberger and P∩olya (1923), P\'olya (1931)) which is very popular in literatures.  The
RRU model is of fundamental importance in many areas of
applications, for instance in economics (c.f., Erev and Roth (1998),
Beggs ( 2005), Hopkins and Posch (2005)), in information science (c.f., Martin and Ho ( 2002)), in resampling theory etc. In clinical trial studies,  the RRU model is  utilized to define a response-adaptive design focusing to   reduce the expected number of patients receiving inferior
treatments  (c.f, Durham, Flournory, Li (1998),   Li, Durham and Flournory (1996), Melifer, Panganoni and
Secchi (2006a,b),  Peganoni and Secchi (2007), May and Flournory (2009) etc).

 Suppose the reinforcement distributions $\mu_1$ and $\mu_2$  have  bounded supports. In   Melifer, Panganoni and Secchi (2006a),   it is showed that the sequence $\{Z_n=Y_{n,1}/|\bm Y_n|\}$ of the random sample proportions in the urn converges to almost surely to a random limit $Z_{\infty}\in [0,1]$. When   $\mu_1= \mu_2$, Crimaldi (2009) proved a central limit theorem by showing almost-sure conditional convergence to a Gaussian kernel of the sequence $\{\sqrt{n}(Z_n-Z_{\infty})\}$. Aletti, May and Secchi (2009) extended Crimaldi's result to a general case that reinforcement means $m_1$ and $m_2$ are equal and proved that $Z_{\infty}$ has no point masses in $[0,1]$ by using this kind of conditional central limit theorem.
 When the means $m_1$ and $m_2$ are different,  the limit proportion $Z_{\infty}$ of a RRU is showed to be a point mass either $1$ and $0$ by Beggs (2005), Hopkins and Posch (2005) and Melifer, Panganoni and Secchi (2006a) under the assumption that the supports of $\mu_1$ and $\mu_2$ are bounded from $0$, and by Aletti, May and Secchi (2009) only under the assumption that  $\mu_1$ and $\mu_2$ have bounded supports. May and Flournory (2009) proved that the sequence $\{Y_{n,1}/Y_{n,2}^{m_1/m_2}\}$ converges to almost surely to a random limit $\psi_{\infty}\in(0,\infty)$ both when $m_1=m_2$ and $m_1\ne m_2$.

 The purpose of this paper is to establish the Gaussian process approximation of the sequence $\{Z_n\}$ when $m_1=m_2$ as well as the sequence $\{Y_{n,1}/Y_{n,2}^{m_1/m_2}\}$ when $m_1\ne m_2$,  under the assumption that $\mu_1$ and $\mu_2$ have only finite $(2+\epsilon)$-th moments. This assumption is much weaker than that  $\mu_1$ and $\mu_2$ have  bounded supports.  We will show that both these sequences can be approximated by a tail stochastic integral with respect to a Brownian motion mixed with a random variable.  It is interesting that, as we will find, the mixed Gaussian process for approximating is nearly independent of the urn composition to be approximated.    Our Gaussian process approximation enables us (i) to establish the law of the iterated logarithm; (ii) to establish the functional limit central limit theorem in both the stable convergence sense and the almost-sure conditional convergence sense; (iii) to prove that the limit $\psi_{\infty}$  (resp. $Z_{\infty}$) has no point masses in $[0,\infty]$ (resp. in $[0,1]$) when $m_1\ne m_2$ (resp. when $m_1= m_2$) under the assumption that $\mu_1$ and $\mu_2$ have only finite $(2+\epsilon)$-th moments. Another implication of our Gaussian approximation is that we are able to establish the central limit theorem in a simple way for the random number $N_{n,k}$ of draws, where $N_{n,k}$ is the number of  type $k$ balls being drawn in the first $n$ samplings. In a response-adaptive design in clinical trials driven by a RRU model, $N_{n,k}$ is the number of patients allocated to treatment $k$, and its asymptotic behaviors are of particular interest.

For the  generalized Friedman urn models, Bai, Hu and Zhang (2002) and Zhang and Hu (2009) established the Gaussian approximation for both the urn proportions  $Y_{n,k}/n$ and the sampling proportion $N_{n,k}/n$.  But the RRU which we
consider here is not covered by their assumptions. The main reason
is that the mean replacement matrix $diag(m_1,m_2)$ of a RRU is not irreducible and hence the limit of $Y_{n,k}/n$ and $N_{n,k}/n$ is not a constant in $(0,1)$.

The paper is organized as follows. The main approximation theorems with applications for equal and unequal reinforcement mean case are stated in Section \ref{section2} and Section \ref{section3}, respectively, and the proofs of the approximations appear in the last section. Some remarks on unsolved problems are discussed in Section \ref{sectionremark}.

   In the sequel of this
paper if having not been specially mentioned, $(U_{l,1}, U_{l,2})$, $l=1,2,\ldots$ are assumed to be independent identically distributed random vectors with finite second moments. Let $X_{m,k}$ be the result of  the $m$-th
drawing, i.e., $X_{m,k}=1$ if the $m$-th drawn ball is of type $k$, and $0$ otherwise. It is obvious that  $ N_{m,k}
=\sum_{j=1}^m  X_{m,k}$ and $X_{m,1}+X_{m,2}=1$.
Denote $\mathcal{F}_n=\sigma\big(U_{l,k},X_{l,k},Y_{l,k}: k=1,2;l=1,\ldots, n\big)$ be the history $\sigma$-field generated by all the observations up to stage $n$, and $\mathcal{F}_{\infty}=\bigvee_n \mathcal{F}_n$.  Further, for two positive sequences $\{a_n\}$ and $\{b_n\}$, we write
$a_n=O(b_n)$ if there is a constant $C$ such that $a_n\le C b_n$,
$a_n\sim b_n$ if $a_n/b_n\to 1$, and $a_n\approx b_n$ if
$a_n=O(b_n)$ and $b_n=O(a_n)$.

%%%%%%%%%%%%%%%%%%%%%%%%%%%%%%%%%%%%%%%%%%%%%%%%%%%%%%%%%%%%%%%
\section{Equal reinforcement mean case }\label{section2}
 \setcounter{equation}{0}

 In this section, we consider the case of $m_1=m_2>0$. Let $\sigma_k^2=\ep[(U_{1,k}/m_k)^2]$, $k=1,2$,
  $$Z_n=\frac{Y_{n,1}}{Y_{n,1}+Y_{n,2}}, \;\; Z_{\infty}(\omega)=\lim Z_n(\omega),\;\; H(\omega)=\frac{\sigma_1^2}{Z_{\infty}}+ \frac{\sigma_2^2}{1-Z_{\infty}}.$$
 To start,  we shall   assume $\pr(Z_{\infty}=0)=\pr(Z_{\infty}=1)=0$, for otherwise $H$ may have no definition.  This result is proved by May and Flournoy (2009) under the condition that the reinforcement distributions $\mu_1$ and $\mu_2$ have bounded supports.  The next theorem tells us that May and Flournoy's condition can be relaxed at least to the assumption of finite second moments.

 \begin{thm}\label{th01} Suppose $m_1=m_2>0$, $\ep U_{1,k}^q<\infty$ for some $q>1$, $k=1,2$. Then the limit $Z_{\infty}$ exists almost surely and $\pr(0<Z_{\infty}<1)=1$.
 \end{thm}

 The following theorem is the main result on the Gaussian approximation.

\begin{thm}\label{th1} Suppose $m_1=m_2>0$, $\ep U_{1,k}^{p}<\infty$, $k=1,2$, where $2\le p<4$.
 Then (possibly in an enlarged probability space) there is standard Brownian motion $B(y)$ such that
  \begin{align}
 Z_{\infty}-Z_n=Z_{\infty}(1-Z_{\infty})H \int_{n H}^{\infty} \frac{ d B(y)}y+ o(\lambda_n)\; a.s.\label{eqth1.1}
 \end{align}
 where
 $$ \lambda_n =\begin{cases} n^{-1/2}(\log\log n)^{1/2}, &\text{ if } p=2\\
 n^{1/p-1}(\log n)^{1/2}, & \text{ if } 2<p<4.
\end{cases}$$
Furthermore, the Brownian motion $B(y)$ can be constructed with a filtration of $\sigma$-fields   $\{\mathscr{G}_n\}$ and a non-decreasing sequence of stopping times $\{T_n\}$ satisfying the following properties:
 \begin{description}
   \item[\rm Property (a)]  $\mathcal{F}_n\subset \mathscr{G}_n$,   $T_n$ is $\mathscr{G}_n$ measurable;
    \item[\rm Property (b)]  $T_n=n H+o(n^{2/p})$ a.s.;
   \item[\rm Property (c)]  Conditional on $\mathcal{F}_n$, $B(T_n+y)-B(T_n)$, $y\ge 0$, is also a standard Brownian motion.
 \end{description}
  \end{thm}

\begin{rem}
 Denote $W(t)= -t\int_{t}^{\infty}y^{-2} dB(y)$. By checking the covariance function, it is easily seen that $W(t), t>0$ is also a standard Brownian motion.
\end{rem}

\begin{rem}
The process in (\ref{eqth1.1}) for approximating is a tail stochastic integral  respective to the Brownian motion. It looks like to be independent of $Z_n$. Actually, according to Property (b) $nH$ can be replaced by $T_n$, and $\sqrt{T_n}\int_{T_n}y^{-1}d B(y)$ is indeed  a normal random variable which is independent of $Z_n$.   We will illustrate this interesting property in Corollary \ref{cor3} in more details.
\end{rem}

We will prove Theorem \ref{th1} by first approximating $Z_n-Z_{\infty}$ to an infinite  summation of a weighted martingale sequence and then approximating the martingale to a Brownian motion by applying the Skorokhod embedding method. The detail proofs of  Theorems \ref{th01} and \ref{th1} will be stated in Section \ref{sectionproofs}. In the sequel of this section, we give several corollaries as applications. Define
$$\widetilde{\sigma}(\omega)=\sqrt{Z_{\infty}(1-Z_{\infty})}\sqrt{(1-Z_{\infty})\sigma_1^2+Z_{\infty}\sigma_2^2},$$
$$\widetilde{\sigma}_n(\omega)=\sqrt{Z_n(1-Z_n)}\sqrt{(1-Z_n)\sigma_1^2+Z_n\sigma_2^2}.$$

The first corollary is  the following law of the iterated logarithm.
  \begin{cor}\label{cor1} Suppose $m_1=m_2>0$, $\ep U_{1,k}^{2}<\infty$, $k=1,2$. Then
  $$ \limsup_{n\to \infty} \frac{ \sqrt{n}(Z_n-Z_{\infty})}{\sqrt{2\log\log n}} =\widetilde{\sigma}\;\; a.s..$$
    \end{cor}

 {\it Proof.} Write $\gamma(x)=\sqrt{x/(2\log\log x)}$, $G(x)=-\int_x^{\infty}y^{-1}d B(y)$. By (\ref{eqth1.1}), we need to show that
 \begin{equation} \label{eqproofLIL1.1}\limsup_{n\to \infty} \gamma(nH)G(nH)=\limsup_{T\to \infty}\gamma(T)|G(T)|=1\; a.s.
  \end{equation}
Note that $\gamma(x)G(x)=xG(x)/\sqrt{2x\log\log x}$, and that $xG(x)$ is also a standard Brownian motion. (\ref{eqproofLIL1.1}) follows from the law of the iterated logarithm of the Brownian motion. $\Box$.

 %%%%%%%%%%%%%%%%%%%%%%%%%%%%%%%%%%%%%%%%%%%%%%%%%%%%%%%%%%%%%%%

 The next corollary is on the functional cental limit theorem.

 \begin{cor}\label{cor2} Suppose $m_1=m_2>0$, $\ep U_{1,k}^{p}<\infty$ for some $p>2$, $k=1,2$.
 Define
 $$ W_n(t)=t\sqrt{n}(Z_{[nt]}-Z_{\infty}), t>0. $$
 Then
 \begin{equation}\label{eqCor2.2.1} W_n(\cdot) \overset{d}\to \widetilde{\sigma} B^{\prime} (\cdot)\;\;\text{stably}, \end{equation}
  in the Skorokhod Topological space $D(0,\infty)$,  where   $B^{\prime}(t)$ is a standard Brownian motion which is independent of $\mathcal{F}_{\infty}$.
  In particular,
 \begin{equation} \label{eqCor2.2.2} \lim_{n\to \infty}\pr\left( \frac{\max_{0\le l\le n}l(Z_l-Z_n)}{\widetilde{\sigma}_n \sqrt{n}}\ge x\right)=e^{-2x^2}, \;\; x>0.
 \end{equation}
Here the stable convergence in  (\ref{eqCor2.2.1}) means  that for any bounded and (uniformly) continues function $f: D(0,\infty)\to (-\infty,\infty)$,
  $$ \ep\left[f\big(W_n(\cdot)\big)I_E\right]\to \ep \left[f\big(\widetilde{\sigma} B^{\prime}(\cdot)\big)I_E\right] \;\; \text{ for any event } E. $$
  \end{cor}

  \begin{rem} The convergece (\ref{eqCor2.2.2})   does not depend on the unknown value of $Z_{\infty}$.
  \end{rem}

{\it Proof.}  (\ref{eqCor2.2.2})  is due to fact that
$$ \max_{0\le l\le n}\frac{l(Z_l-Z_n)}{\widetilde{\sigma} \sqrt{n}}=\sup_{0<t\le 1}\frac{W_n(t)-tW_n(1)}{\widetilde{\sigma}}
\overset{d}\to \sup_{0<t\le 1} (B^{\prime}(t)-tB^{\prime}(1)) $$
and that $B^{\prime}(t)-tB^{\prime}(1)$ is a Brownian bridge. For (\ref{eqCor2.2.1}), note that  $W(x)=-t\int_t^{\infty} y^{-1} d B(y)$ is also  a standard Brownian motion. By (\ref{eqth1.1}),
$$  n(Z_n-Z_{\infty})-\frac{\widetilde{\sigma}}{\sqrt{H}} W(nH)=o(n^{1/2}
)\;\; a.s., $$
which implies that for any $T>0$,
$$\sup_{0<t\le T}\left|W_n(t)-\widetilde{\sigma} \frac{W (nHt)}{\sqrt{nH}}\right|\to o(1)\;a.s..$$
For the Brownian motion $W(\cdot)$, we have
$$  W(n\cdot)/\sqrt{n}\overset{d}\to  B^{\prime}(\cdot)\;\; \text{mixing},$$
i.e.,  for any given event $E$ with $\pr(E)>0$, the conditional distribution of $W(n\cdot)/\sqrt{n}$ converges to a Brownian motion. It follows that
  $$\big(\widetilde{\sigma}/\sqrt{H}, H, W(n\cdot)/\sqrt{n}\big)\overset{d}\to \big(\widetilde{\sigma}/\sqrt{H}, H, B^{\prime}(\cdot)\big)\;\;\text{stably}.$$
  Note that $\widetilde{\sigma} \frac{W(nH\cdot)}{\sqrt{nH}}$ is a continuous function of $(\widetilde{\sigma}/\sqrt{H}, H, W(n\cdot)/\sqrt{n})$ of the form $f(r,h,x(\cdot))=r x(\cdot h)$. It follows that
  $$\widetilde{\sigma} \frac{W(nH\cdot)}{\sqrt{nH}}\overset{d}\to \widetilde{\sigma} \frac{B^{\prime}(\cdot H)}{\sqrt{H}}\overset{d}=\widetilde{\sigma}  B^{\prime}(\cdot)\;\;\text{stably}. $$
  The proof is now completed. $\Box$

\bigskip
Corollary \ref{cor2} implies the central limit theorem for $\sqrt{n}(Z_n-Z_{\infty})$. Aletti, May and Secchi (2009) proved a strong version of the central limit theorem.  For every Borel set $B$, every $\omega$, and $n=1,2,\ldots,$, define
  $$K_n(\omega, B)=\pr\big(\sqrt{n}(Z_n-Z_{\infty})\in B\big|\mathcal{F}_n\big)(\omega), $$
 i.e., $K_n$ is a version of the condition distribution of $\sqrt{n}(Z_n-Z_{\infty})$ given $\mathcal{F}_n$. Aletti, May and Secchi (2009) showed that,  if $m_1=m_2$ and the distributions of $\mu_1$ and $\mu_2$ have bounded supports, then for almost every $\omega$, the sequence of probability distributions $K_n(\omega,\cdot)$ converges weakly to the normal distribution
   $$ N\big(0, \widetilde{\sigma}^2(\omega)\big). $$
   We denote this kind of convergence by
  \begin{equation}\label{eqCor2.3.1}
   \sqrt{n}(Z_n-Z_{\infty})\Big|_{\mathcal{F}_n}\overset{d}\to N\big(0, \widetilde{\sigma}^2(\omega)\big)\; a.s..
   \end{equation}
This kind of conditional central limit theorem was first established by Grimadi (2008) for a special case that $\mu_1=\mu_2$. Aletti, May and Secchi (2009) also showed that (\ref{eqCor2.3.1}) implies that $Z_{\infty}$ has no point masses in $(0,1)$.  Our next corollary tells that (\ref{eqCor2.3.1}) and a type of conditional functional central limit theorem are  conclusions of the Gaussian approximation.

  \begin{cor}\label{cor3} Suppose $m_1=m_2>0$, $\ep U_{1,k}^{p}<\infty$ for some $p>2$, $k=1,2$. Then
  \begin{equation}\label{eqCor2.3.2}
  \sup_{t\ge 1}\left| \sqrt{n}(Z_{\infty}-Z_{[nt]})-\widetilde{\sigma}_n\sqrt{T_n}\int_{T_nt}^{\infty}\frac{d B(y)}{y}\right|
  =o(n^{-\epsilon})\;\; a.s.
  \end{equation}
  for some $\epsilon>0$,    and further, $\overline{B}_n(t)=-t\sqrt{T_n}\int_{T_nt}^{\infty} y^{-1} dB(y)$, $t\ge 1$, is also a standard Brownian on $[1,\infty)$ which is independent of $\mathcal{F}_n$.

  As a consequence,
 \begin{equation}\label{eqCor2.3.3} W_n(\cdot)\big|_{\mathcal{F}_n} \overset{d}\to \widetilde{\sigma}(\omega) B^{\prime} (\cdot)\;\; a.s.,
 \end{equation}
  in the Skorokhod Topological space $D[1,\infty)$,  where   $W_n(t)$ is defined as in Corollary \ref{cor2}, $B^{\prime}(t)$ is a standard Brownian motion which is independent of $\mathcal{F}_{\infty}$.

  In particular, (\ref{eqCor2.3.1}) holds,   $Z_{\infty}$ has no point masses in $(0,1)$, and there exists a sequence of standard normal random variables for which $\zeta_n$   is independent of $\mathcal{F}_n$  and
  \begin{equation}\label{eqCor2.3.4}
  \sqrt{n}(Z_n-Z_{\infty})=\widetilde{\sigma}_n \zeta_n+o(n^{-\epsilon})\;\; a.s.,
  \end{equation}
    \end{cor}

  {\it Proof.} We first prove (\ref{eqCor2.3.2}). By (\ref{eqth1.1}),
 $$ \sup_{t\ge 1}\left|\sqrt{n}(Z_{\infty}-Z_{[nt]})-Z_{\infty}(1-Z_{\infty})H\sqrt{n} \int_{nt H}^{\infty} \frac{ d B(y)}y\right|= o(n^{-\epsilon})\;a.s. $$
  Let $H_n=\sigma_1^2/Z_n + \sigma_2^2/(1-Z_n).$
Note that
$Z_n-Z_{\infty}=O(\sqrt{n^{-1}\log\log n})$ a.s. by   Corollary \ref{cor1}, and $T_n/n=H+o(n^{2/p-1})$ a.s..  It follows that
$\sqrt{n}Z_{\infty}(1-Z_{\infty})H-\widetilde{\sigma}_n\sqrt{T_n}=o(n^{2/p-1/2})$ a.s.. Further,
$$ n^{2/p-1/2}\sup_{t\ge 1}\left|  \int_{nt H}^{\infty} \frac{ d B(y)}y\right|=n^{2/p-1/2}O(\sqrt{n^{-1}\log\log n})=o(n^{-\epsilon}) \;\;a.s. $$
It remains to show that
$$\sup_{t\ge 1}\left|  \int_{nt H}^{T_nt} \frac{ d B(y)}y\right|=o(n^{-1/2-\epsilon})\;\; a.s. $$
 Write  $a_n=T_n-nH$. Note that $W(x)=-x\int_x^{\infty}y^{-1}d B(y)$ is a standard Brownian motion and
 $$\int_{nt H}^{T_nt} \frac{ d B(y)}y=\frac{W(nHt+a_nt)-W(nHt)}{T_nt}+
 \frac{W(ntH)}{ntH}\frac{nH-T_n}{T_n}. $$
The second term on the right hand of the above equality does not exceed
$$O(\sqrt{n^{-1}\log\log n})o(n^{2/p-1}) $$
 uniformly in $t\ge 1$ almost surely, by the law of the iterated logarithm. The first term does not exceed
$$ O\left(\frac{\sqrt{a_nt(\log(nH)+\log\log (nHt))}}{T_n t}\right)=o(n^{1/p-1}(\log n)^{1/2})  $$
uniformly in $t\ge 1$ almost surely, by the path properties of a Brownian motion (c.f. Hanson and Russo (1983)). The proof of (\ref{eqCor2.3.2}) is now proved.

Let $B_n(y)=B(T_n+y)-B(T_n)$.  Then conditional on $\mathscr{G}_n$, $B_n(y)$ is a standard Brownian motion. It is obvious that
$$   \overline{B}_n(t) =-t\sqrt{T_n}\int_{T_n(t-1)}^{\infty}\frac{d B_n(y)}{T_n+y}. $$
 Hence, conditional on $\mathscr{G}_n$, $\overline{B}_n(t)$, $t\ge 1$,  is a mean zero Gaussian process with covariance function
 $$ T_n ts \int_{T_n(t-1)}^{\infty} \frac{d y}{(T_n+y)^2}=s \;\; \text{ for } t\ge s\ge 1. $$
It follows that $\overline{B}_n(t)$, $t\ge 1$, is a standard Brownian motion and is independent of $\mathscr{G}_n$. So, it is independent of $\mathcal{F}_n$ because $\mathcal{F}_n\subset \mathscr{G}_n$. The proof of the main part of the corollary is now completed.

Now for (\ref{eqCor2.3.3}), from (\ref{eqCor2.3.2}) it follows that
$$ dist\big(W_n(\cdot),\widetilde{\sigma}_n B_n(\cdot)\big)\to 0 \;\;a.s.\; \text{ in } D[1,\infty), $$
where $dist(\cdot,\cdot)$ is a metric in $D[1,\infty)$.
Note that $\widetilde{\sigma}_n$ is $\mathscr{F}_n$-measurable and $\widetilde{\sigma}_n\to\widetilde{\sigma}$ a.s., and, conditional on $\mathcal{F}_n$,   $\widetilde{\sigma}_nB_n(\cdot)$ and $\widetilde{\sigma}_n B^{\prime}(\cdot)$ has the same distribution. Hence, (\ref{eqCor2.3.3}) follows from (\ref{eqCor2.3.4}) by  noting  the following fact that:
  \begin{align*}
  &\xi_n(\cdot)\big|_{\mathcal{F}_n}\overset{d}\to  \widetilde{\sigma}B^{\prime}(\cdot)\; a.s.
   \text{ and } dist\big(\xi_n(\cdot),\eta_n(\cdot)\big)\to 0\; a.s. \;\; \text{ in } D[1,\infty)\\
 &\quad  \Longrightarrow
   \eta_n(\cdot)\big|_{\mathcal{F}_n}\overset{d}\to \widetilde{\sigma}B^{\prime}(\cdot)\; a.s.\;\; \text{ in } D[1,\infty).
  \end{align*}
  This  fact follows from that for any bounded and uniformly continuous function $f:D[1,\infty)\to (-\infty,\infty)$,
  \begin{align*}
  & \limsup_{n\to\infty}\Big|\ep\big[f(\eta_n)-f(\xi_n)\big|\mathcal{F}_n\big]\Big|\\
  \le &
 \ep\Big[\limsup_{n\to\infty}|f(\eta_n)-f(\xi_n)|\Big|\bigvee_n\mathcal{F}_n\Big]=0 \;\; a.s.,
 \end{align*}
 due to Lemma A.2 of Crimaldi (2009).

Finally, (\ref{eqCor2.3.4}) follows from (\ref{eqCor2.3.2}) by letting $\zeta_n=\overline{B}_n(1)$, and (\ref{eqCor2.3.1}) is a conclusion of (\ref{eqCor2.3.4}) or (\ref{eqCor2.3.3}).  Aletti, May and Secchi (2009)  showed that (\ref{eqCor2.3.1}) implies  $Z_{\infty}$ having no point masses in $(0,1)$
by utilizing  a metric of the weak convergence of probability measures with  the limit distribution being absolutely continuous. Here we give a straightforward proof. Let $f(t)=e^{-t^2/2}$  be the characteristic function of a standard normal distribution.   Firstly, note that   (\ref{eqCor2.3.1}) implies that for every $\bigvee_n \mathcal{F}_n$-measurable event $E$,
$$ \lim_{n\to\infty}\ep[e^{it\sqrt{n}(Z_n-Z_{\infty})}I_E|\mathcal{F}_n]=f(\widetilde{\sigma}t)I_E\;\; a.s. $$
In fact, if let $I_n=\ep[I_E|\mathcal{F}_n]$, then $I_n\to I_E$ a.s.. And hence
\begin{align*}
&\lim_{n\to\infty}\ep[e^{it\sqrt{n}(Z_n-Z_{\infty})}I_E|\mathcal{F}_n]=
\lim_{n\to\infty}\ep[e^{it\sqrt{n}(Z_n-Z_{\infty})}I_n|\mathcal{F}_n]\\
&=\lim_{n\to\infty}\ep[e^{it\sqrt{n}(Z_n-Z_{\infty})}|\mathcal{F}_n]I_n
=f(\widetilde{\sigma}t)I_E\;\; a.s.,
\end{align*}
where in the fist equality we use the fact that
$$\eta_n\to 0 \;\;a.s. \; \text{ and } |\eta_n|\le M\; a.s. \Longrightarrow \ep[\eta_n|\mathcal{F}_n]\to 0\;\; a.s.. $$  This fact is due to Lemma A.2 of Crimaldi (2009). Next, choosing  $E=\{Z_{\infty}=p\}$,  $p\in (0,1)$, yields
\begin{align*}
&f(\widetilde{\sigma}t)I_E= \lim_{n\to\infty}\ep[e^{it\sqrt{n}(Z_n-p)}I_E|\mathcal{F}_n]\\
=&\lim_{n\to\infty}e^{it\sqrt{n}(Z_n-p)}\ep[I_E|\mathcal{F}_n]
=\lim_{n\to\infty}e^{it\sqrt{n}(Z_n-p)}I_E\;\; a.s..
\end{align*}
Hence, $|f(\widetilde{\sigma}t)|I_E=I_E$ a.s.. So $I_E=0$ a.s. because $|f(\widetilde{\sigma}t)|<1$ on $E$.
The proof is now completed. $\Box$.

From the above proof, we obtain the following corollary.
\begin{cor} Suppose $\{Y_\infty, Y_n,n\ge 1\}$ is a sequence of random variables, $\{a_n\}$ is a sequence of constants with $a_n\to \infty$, and $\mathscr{G}_n$ is a filtration of $\sigma$-fields such that $Y_n$ is $\mathscr{G}_n$-measurable. If for almost every $\omega$, the distribution of $a_n(Y_n-Y_{\infty})$ conditional on $\mathscr{G}_n$ converges to a non-degenerate distribution, then $Y_{\infty}$ has no point masses.
\end{cor}

\begin{rem}By (\ref{eqCor2.3.1}), it can also be shown   that for every event $E$,
\begin{align*}
&\limsup_{|t|\to\infty}\limsup_{n\to\infty}\left|\ep\Big[\exp\big\{-itn^{1/2}Z_{\infty}\big\}I_E\big]\right|\\
\le & \limsup_{|t|\to\infty}\limsup_{n\to\infty}\ep\left| \ep\Big[\exp\big\{itn^{1/2}(Z_n-Z_{\infty})\big\}I_E\big|\mathcal{F}_n\big] \right|\\
=&\limsup_{|t|\to\infty}\limsup_{n\to\infty}\ep\left|\ep\big[f\big(t \widetilde{\sigma}_n\big)I_E\big|\mathcal{F}_n\big]\right|\\
=&\limsup_{|t|\to\infty} \ep\big[\left|f\big(t \widetilde{\sigma}_{\infty}\big)\right|I_E \big]=0.
\end{align*}
If we denote $f_E(t)$ be the characteristic function of the conditional distribution of $Z_{\infty}$ given $E$ with $\pr(E)>0$. Then the above equality means that
$$ \lim_{|t|\to \infty} \limsup_{n\to \infty}\left|f_E(tn^{1/2})\right|=0. $$
Note that for any $t_0\ge 1$, $n_0\ge 1$ and $|s|\ge 2t_0n_0$, there exist a real number $t$ with $|t|\ge t_0$  and an integer $n\ge n_0$ such that $s=tn^{1/2}$. We conclude that
$$\lim_{|s|\to \infty} f_E(s)=0. $$
This is related to the Cram\'er condition. Obviously, if $E=\{Z_{\infty}=p\}$, $p\in (0,1)$, and $\pr(E)>0$, then
$|f_E(t)|\equiv 1$ which is a contradiction.
\end{rem}

\bigskip

%%%%%%%%%%%%%%%%%%%%%%%%%%%
The next corollary is the central limit theorem for the random number of draws.

\begin{cor} \label{Cor4} Suppose $m_1=m_2>0$, $\ep U_{1,k}^{p}<\infty$ for some $p>2$, $k=1,2$.  Then
\begin{equation}
\sqrt{n}\left(\frac{N_{n,1}}{n}-Z_{\infty}\right)\overset{d}\to h(\omega) \cdot N(0,1) \;\; \text{stably},
\end{equation}
where $h(\omega)= \sqrt{Z_{\infty}(1-Z_{\infty})}\sqrt{(1-Z_{\infty})(2\sigma_1^2-1)+Z_{\infty}(2\sigma_2^2-1)}$, and $N(0,1)$ is a standard normal random variable which is independent of $\mathcal{F}_{\infty}$.
\end{cor}

{\it Proof.} We need to prove
\begin{align*}
&\sqrt{n}\left(\frac{N_{n,1}}{n}-Z_{\infty}\right)\\
\overset{d}\longrightarrow &-\sqrt{Z_{\infty}}(1-Z_{\infty})N_1\Big(0,\sigma_1^2-1\Big)\\
&\quad +\sqrt{1-Z_{\infty}}Z_{\infty}N_2\Big(0,\sigma_2^2-1)\Big)+\widetilde{\sigma} \cdot N_3(0,1)
\;\; \text{stably},
\end{align*}
where $N_k\Big(0,\sigma_k^2-1\Big)$, $k=1,2$, $N_3(0,1)$ are three independent normal random variables which are independent of $\mathcal{F}_{\infty}$. Write
\begin{equation}\label{eqproofCor2.4.1} A_{n,k}=\frac{\sum_{l=1}^nX_{l,k}(U_{l,k}/m_k-1)}{N_{n,k}}, \;\; k=1,2.
\end{equation}
Then
$$ A_{n,k}=O\big(\sqrt{n^{-1}\log\log n}\big)\; a.s.\;\text{ and }\;\; \frac{Y_{n,k}}{m_kN_{n,k}}=1+A_{n,k}. $$
By the Taylor expansion and (\ref{eqCor2.3.4}), we have that
\begin{align}\label{eqproofCor2.4.2}
&\frac{N_{n,1}}{n}-Z_{\infty}=\frac{N_{n,1}}{N_{n,1}+N_{n,2}}-Z_{\infty}
=\frac{\frac{Z_n}{1+A_{n,1}}}{\frac{Z_n}{1+A_{n,1}}+\frac{1-Z_n}{1+A_{n,2}}}-Z_{\infty}\nonumber\\
=& -Z_n(1-Z_n)\big(A_{n,1}-A_{n,2})+(Z_n-Z_{\infty})+O\Big(\frac{\log\log n}{n}\Big)\nonumber\\
=&-\frac{Z_n(1-Z_n)}{\sqrt{N_{n,1}}}\big(\sqrt{N_{n,1}}A_{n,1}\big)
+\frac{Z_n(1-Z_n)}{\sqrt{N_{n,2}}}\big(\sqrt{N_{n,2}}A_{n,2}\big)\nonumber\\
&+\widetilde{\sigma}_n\zeta_n+o(n^{-1/2})\; a.s..
\end{align}
  Note that $\zeta_n$ is a standard normal random variable which is independent of $\widetilde{\sigma}_n$, $Z_n$, $N_{n,k}$ and $A_{n,k}$, $k=1,2$. Also,  $\widetilde{\sigma}_n\to \widetilde{\sigma}$ a.s., $Z_n(1-Z_n)\sqrt{n}/\sqrt{N_{n,1}}\to \sqrt{Z_{\infty}}(1-Z_{\infty})$ a.s. and $Z_n(1-Z_n)\sqrt{n}/\sqrt{N_{n,2}}\to Z_{\infty}\sqrt{1-Z_{\infty}}$. The proof is completed if we have shown that
\begin{equation}\label{eqproofCor2.4.3}  \big(\sqrt{N_{n,1}}A_{n,1},  \sqrt{N_{n,2}}A_{n,2}\big)\overset{d} \longrightarrow
\Big(N_1\big(0,\sigma_1^2-1\big), N_2\big(0,\sigma_2^2-1\big)\Big)\;\; \text{mixing}.
\end{equation}
Note that $\sigma_k^2-1=\Var\big(U_{1,k}/m_k\big)$, $k=1,2$. The above convergence follows from Theorem 4.1 of May and Flournoy (2009).  $\Box$

The next corollary tells us that conditional on $\mathcal{F}_n$, the conditional distribution of $\sqrt{n}(N_{n,1}/n-Z_{\infty})$ does not converge.
\begin{cor} \label{Cor5} Suppose $m_1=m_2>0$, $\ep U_{1,k}^{p}<\infty$ for some $p>2$, $k=1,2$. Let $E$ be an event that for $\omega\in E$ there is a distribution $F_{\omega}$ for which
\begin{equation}\label{eqCor2.5.1}
\sqrt{n}\big(\frac{N_{n,1}}{n}-Z_{\infty}\big)\Big|_{\mathcal{F}_n}\overset{d}\to F_{\omega}.
\end{equation}
Then $\pr(E)=0$.
\end{cor}

{\it Proof.} Recall (\ref{eqproofCor2.4.2}). Let $\eta_n=\sqrt{nZ_n(1-Z_n)}(A_{n,1}-A_{n,2})$. Note that $Z_n\to Z_{\infty}$ a.s.,
$$ \widetilde{\sigma}_n\xi_n\big|_{\mathcal{F}_n}\overset{d}\to \widetilde{\sigma} N(0,1)\;\; a.s. $$
and that $Z_n$, $A_{n,1}$ and $A_{n,2}$ are $\mathcal{F}_n$-measurable. By (\ref{eqproofCor2.4.2}) and (\ref{eqCor2.5.1}), there exists an event $\Omega_0$ with $\pr(\Omega_0)=1$ such that $\sqrt{n}Z_n(1-Z_n)(A_{n,1}-A_{n,2})$ converges on $E\cap \Omega_0$.  So there exists a random variable $\eta$ such that
$$ \eta_n(\omega)\to \eta(\omega)\;\; \forall \omega\in E\cap \Omega_0. $$
Suppose $\pr(E)>0$. Choose $x$ such that $\pr(\eta>x, E)>0$. Then it follows that
$$\pr(\eta_n\le x,\eta>x,E)\to \pr(\eta\le x,\eta>x,E)=0.$$
So
$ \pr(\xi_n\le x|\xi>x,E)\to 0. $
On the other hand, according to (\ref{eqproofCor2.4.3}) we have
$$ \eta_n=\frac{\sqrt{nZ_n(1-Z_n)}}{\sqrt{N_{n,1}}}\sqrt{N_{n,1}}A_{n,1}-\frac{\sqrt{nZ_n(1-Z_n)}}{\sqrt{N_{n,2}}}\sqrt{N_{n,2}}A_{n,2}
\overset{d}\to N(0,1) \;\;\text{mixing}, $$
where $N(0,1)$ is independent of $\mathcal{F}_{\infty}$.
It follows that
$$ \lim_n\pr(\eta_n\le x|\eta>x,E) =\Phi(x)>0. $$
We get a contradiction. The proof is completed. $\Box$

%%%%%%%%%%%%%%%%%%%%%%%%%%%%%%%%%%%%%%%%%%%%%%%%%%%%%%%%%%%%%%%%%%%%%%%%%%%%%%%%%%%%%%%%%5

\section{Unequal reinforcement mean case }\label{section3}

In this section, we consider the case of $m_1\ne m_2$.  Without loss of generality, we assume that $0<m_1< m_2$. Denote
$\rho=m_1/m_2$, and
$$ \psi_n= \frac{Y_{n,1}}{Y_{n,2}^{\rho}},\;\; \psi_{\infty}=\lim_{n\to\infty} \psi_n. $$
May and Flournoy (2009) proved that the limit $\psi_{\infty}$ exists almost surely with $\pr(0<\psi_{\infty}<\infty)=1$ when the reinforcement distributions $\mu_1$ and $\mu_2$ have bounded supports. Durham and Yu (1990) proved a similar result as that
$$ \frac{N_{n,1}}{N_{n,2}^{\rho}} \text{ converges almost surely to a finite limit } \eta_{\infty}. $$
It is easily seen that
$$ \eta_{\infty}=\frac{m_2^{\rho}}{m_1}\psi_{\infty}\;\; a.s. $$
and
$$ \lim_{n\to \infty}\frac{N_{n,1}}{n^{\rho}}=\frac{m_2^{\rho}}{m_1}\psi_{\infty}\; a.s.,\;\;
\lim_{n\to \infty}\frac{Y_{n,1}}{n^{\rho}}= m_2^{\rho} \psi_{\infty}\; a.s. $$
In a recent manuscript of Zhang et al (2010), it is  proved that the weakest condition for $\pr(0<\psi_{\infty}<\infty)=1$ is that $\ep[U_{1,k}\log^+U_{1,K}]<\infty$, $k=1,2$, and a general multi-color RRU is consider.  For the completeness of this paper, we will give a simple proof under the assumption of finite $(1+\epsilon)$-th moments for the two-color case. The following is the result.

\begin{thm}\label{th02} Suppose $\ep U_{1,k}^q<\infty$ for some $q>1$, $m_k>0$,  $k=1,2$. Then the limit $\psi_{\infty}$ exists almost surely  and $\pr(0<\psi_{\infty}<\infty)=1$  both when  $m_1=m_2$ and $m_1\ne m_2$.
\end{thm}

The following theorem is our main result on the Gaussian process approximation for $\psi_n$. From the Gaussian approximation we are able to show that $\psi_{\infty}$ has no point masses in $(0,\infty)$. And accordingly, all the limits of the sequences
$\{Y_{n,1}/Y_{n,2}^{\rho}\}$, $\{Y_{n,1}/n^{\rho}\}$, $\{Y_{n,1}/|\bm Y_n|^{\rho}\}$, $\{N_{n,1}/N_{n,2}^{\rho}\}$  and $\{N_{n,1}/n^{\rho}\}$ have no point masses in $[0,\infty]$.

\begin{thm}\label{th2} Suppose $m_2>m_1>0$,  $\ep U_{1,k}^{p}<\infty$, $k=1,2$, for some $p> 2$. Denote $\sigma_k^2=\ep[(U_{1,k}/m_k)^2]$, $k=1,2$. Let $\delta_0=\min\{(1-\rho)/\rho,1/2-1/p\}$.
 Then (possibly in an enlarged probability space) there is standard Brownian motion $B(y)$ such that for any $0<\delta<\delta_0$,
  \begin{align}
  \psi_{\infty}-\psi_n = &\frac{\sigma_1\sqrt{m_1  }}{m_2^{\rho/2}}
     \int_{n^{\rho}/\psi_{\infty}}^{\infty}  \frac{ d B(y)}{y}+ o( n^{ -\rho(1+\delta)/2})\label{eqth2.1}\\
=& -n^{-\rho/2}\frac{\sigma_1\sqrt{m_1 \psi_{\infty}}}{m_2^{\rho/2}}
   \frac{W\left(n^{\rho}/\psi_{\infty}\right)}{\sqrt{n^{\rho}/\psi_{\infty}}}+ o( n^{-\rho(1+\delta)/2})\; a.s.,\label{eqth2.2}
 \end{align}
 where $W(x)= -x\int_{x}^{\infty}y^{-1} dB(y)$ is  also a standard Brownian motion.

 Furthermore, the Brownian motion $B(y)$ can be constructed with a filtration of $\sigma$-fields   $\{\mathscr{G}_n\}$ and a non-decreasing sequence of  stopping times $\{T_n\}$ satisfying Properties (a) and  (c) in Theorem \ref{th1}, and
 \begin{description}
   \item[\rm Property (b$^{\prime}$)]  $T_n=n^{\rho}/\psi_{\infty} +o(n^{\rho(1-\delta)})$ a.s. for  $0<\delta<\delta_0$.
  \end{description}
     \end{thm}

The proof of this theorem will be given in the last section. Next, we state  several corollaries. The first one is on the law of iterated logarithm and the central limit theorem for $\psi_n$.

  \begin{cor}\label{cor3.1}
  Under the conditions in Theorem \ref{th2},
  \begin{equation}\label{eqCor3.1.1}
  \limsup_{n\to \infty}\frac{n^{\rho/2}\big(\psi_n-\psi_{\infty}\big)}{\sqrt{2\log\log n}} = \frac{\sigma_1\sqrt{m_1 \psi_{\infty}}}{m_2^{\rho/2}}\;\; a.s.
   \end{equation}
  and there exists a sequence $\{\zeta_n\}$ of standard normal random variables for which $\zeta_n$ is independent of $\mathcal{F}_n$ and
  \begin{equation}\label{eqCor3.1.2}
  n^{\rho/2}\big(\psi_n-\psi_{\infty}\big)=\frac{\sigma_1\sqrt{m_1 \psi_n}}{m_2^{\rho/2}} \zeta_n+o(n^{-\epsilon}) \; a.s. \text{ for some  } \epsilon>0.
  \end{equation}
  Hence
  \begin{equation}\label{eqCor3.1.3}n^{\rho/2}\big(\psi_n-\psi_{\infty}\big)\Big|_{\mathcal{F}_n}\overset{d}\longrightarrow N\Big(0,\frac{\sigma_1^2m_1}{m_2^{\rho}}\psi_{\infty}(\omega)\Big)\;\; a.s.
  \end{equation}
  and $\psi_{\infty}$ has no point masses in $(0,\infty)$.
  \end{cor}

 {\it Proof.}  (\ref{eqCor3.1.1}) follows from (\ref{eqth2.2}) and the law of iterated logarithm of the Brownian motion. (\ref{eqCor3.1.2}) can be proved in the same way as proving  Corollary \ref{cor1}. $\Box$

\begin{cor}\label{cor3.2}
  Under the conditions in Theorem \ref{th2},
    \begin{equation}\label{eqCor3.2.1}n^{\rho/2}\big(\frac{N_{n,1}}{N_{n,2}^{\rho}} -\eta_{\infty}\big)\overset{d}\longrightarrow N (0,1 )\cdot\sqrt{\eta_{\infty}(2\sigma_1^2-1)}\;\; \text{ stably}
  \end{equation}
and
    \begin{equation}\label{eqCor3.2.2}
        \begin{cases} n^{\rho/2}\big(\frac{N_{n,1}}{n^{\rho}} -\eta_{\infty}\big)\overset{d}\to  N (0,1 )\cdot\sqrt{\eta_{\infty}(2\sigma_1^2-1)}\;\; \text{ stably} & \text{ if } \rho<2/3\\
     n^{1-\rho}\big(\frac{N_{n,1}}{n^{\rho}} -\eta_{\infty}\big)\to -\rho \eta_{\infty}^2 \;\; a.s.  & \text{ if } \rho>2/3,
    \end{cases}
  \end{equation}
  \begin{equation}\label{eqCor3.2.3}
 n^{1-\rho}\big(1-\frac{N_{n,2}}{n} \big)\to   \eta_{\infty} \;\; a.s.
  \end{equation}
     \end{cor}

{\it Proof.}
   For (\ref{eqCor3.2.1}), let $A_{n,k}$ be defined as in (\ref{eqproofCor2.4.1}). Then
 $$A_{n,2}=O(\sqrt{N_{n,2}^{-1}\log\log N_{n,2}})=O(\sqrt{n^{-1}\log\log n})=o(n^{-\rho/2-\epsilon})\;\;a.s., $$
 $$A_{n,1}=O(\sqrt{N_{n,1}^{-1}\log\log N_{n,1}})=O(\sqrt{n^{-\rho}\log\log n})\;\;a.s.. $$
 Note that $\eta_{\infty}=\psi_{\infty} m_2^{\rho}/m_1$. It follows that
 \begin{align*}
  &\frac{N_{n,1}}{N_{n,2}^{\rho}}-\eta_{\infty}= -\eta_{\infty}+\psi_n \frac{m_2^{\rho}}{m_1}\frac{(1+A_{n,2})^{\rho}}{1+A_{n,1}}\\
  =&-\eta_{\infty}+\psi_n \frac{m_2^{\rho}}{m_1}(1-A_{n,1})+o(n^{-\rho/2-\epsilon})\\
 = & (\psi_n-\psi_{\infty})\frac{m_2^{\rho}}{m_1}-n^{-\rho/2}\psi_n \frac{m_2^{\rho}}{m_1}\sqrt{n^{\rho}/N_{n,1}}  \big(\sqrt{N_{n,1}}A_{n,1}\big)+o(n^{-\rho/2-\epsilon})\\
 =&n^{-\rho/2}\Big\{\sqrt{\psi_n m_1/m_2^{\rho}}\; \zeta_n-\psi_n \frac{m_2^{\rho}}{m_1}\sqrt{n^{\rho}/N_{n,1}}  \big(\sqrt{N_{n,1}}A_{n,1}\big)+o(n^{-\epsilon})\Big\}\;\; a.s.
 \end{align*}
 The proof of (\ref{eqCor3.2.1}) is completed by noting that $\psi_n m_1/m_2^{\rho}\to \eta_{\infty}$ a.s., $N_{n,1}/n^{\rho}\to \eta_{\infty}$ a.s.,
$\sqrt{N_{n,1}}A_{n,1}\overset{d}\to N(0,\sigma_1^2-1)$ mixing, and $\zeta_n$ is a standard normal random variable which is independent of $\mathcal{F}_n$.

 For (\ref{eqCor3.2.2}), it is sufficient to note that
 \begin{align*}
 \frac{N_{n,1}}{n^{\rho}}=&\frac{N_{n,1}}{N_{n,2}^{\rho}}\left(1-\frac{N_{n,1}}{n}\right)^{\rho}
 =\frac{N_{n,1}}{N_{n,2}^{\rho}}\left(1-\rho\frac{N_{n,1}}{n}+O\Big(\frac{N_{n,1}}{n}\Big)^2\right)\\
 =&\frac{N_{n,1}}{N_{n,2}^{\rho}}-\rho\frac{N_{n,1}}{N_{n,2}^{\rho}}\frac{N_{n,1}}{n^{\rho}}n^{\rho-1}+O(n^{2(\rho-1)})\\
 =&\frac{N_{n,1}}{N_{n,2}^{\rho}}-\rho\eta_{\infty}^2n^{\rho-1}+o(n^{\rho-1})\;\; a.s.
 \end{align*}
 (\ref{eqCor3.2.3}) is obvious because
 $$1-\frac{N_{n,2}}{n}=\frac{N_{n,1}}{n}\sim \frac{\eta_{\infty}n^{\rho}}{n}\;\; a.s. $$
 The proof is now completed. $\Box$.

Finally, we give the functional central limit theorem.
 \begin{cor}\label{cor3.3}
 Define
 $$ W_n(t)=n^{\rho/2}t^{\rho}(\psi_{[nt]}-\psi_{\infty}), t>0. $$
 Then
 \begin{equation}\label{eqCor3.3.1} W_n(t) \overset{d}\to \sigma_1\sqrt{\eta_{\infty}} B^{\prime} (t^{\rho})\;\;\text{stably}, \end{equation}
  in the Skorokhod Topological space $D(0,\infty)$,  where   $B^{\prime}(t)$ is a standard Brownian motion which is independent of $\mathcal{F}_{\infty}$.
  In particular,
 \begin{equation} \label{eqCor3.3.2} \lim_{n\to \infty}\pr\left( \frac{\max_{0\le l\le n}l^{\rho}(\psi_l-\psi_n)}{\sigma_1\sqrt{N_{n,1}}}\ge x\right)=e^{-2x^2}, \;\; x>0.
 \end{equation}
  \end{cor}

{\it Proof.} The proof of (\ref{eqCor3.3.1}) is similar to that of (\ref{eqCor2.2.1}) by noting that $\eta_{\infty}=\psi_{\infty}m_1/m_2^{\rho}$. For (\ref{eqCor3.3.2}), it is sufficient to see that
$$ \max_{0\le l\le n}\frac{l^{\rho}(\psi_{\infty}-\psi_n)}{n^{\rho/2}}=\sup_{0<t\le 1}(W_n(t)-t^{\rho}W_n(1))
 $$
and $N_{n,1}\sim \eta_{\infty}n^{\rho}$ a.s.. $\Box$

 %%%%%%%%%%%%%%%%%%%%%%%%%%%%%%%%%%%%%%%%%%%%%%%%%%%%%%%%%%

\section{Concluding Remark}\label{sectionremark}

We approximated $Z_n-Z_{\infty}$ and $\psi_n-\psi_{\infty}$ by a kind of Gaussian process $\int_t^{\infty}y^{-1} d B(y)$, which is a tail stochastic integral with respective to a Brownian motion,   with time $t$ stopping at a random variable $nH$ or $n/\psi_{\infty}$, where $H^2=(1+\psi_{\infty})(\sigma_1^2/\psi_{\infty}+\sigma_2^2)$. But this does not mean that $\int_{nH}^{\infty}y^{-1} d B(y)$ and $\int_{n/\psi_{\infty}}^{\infty}y^{-1} d B(y)$ are     Gaussian random variables and their distributions are unknown because the mixing distribution of $\psi_{\infty}$ is unknown. For deriving the asymptotic distributions, the approximations (\ref{eqCor2.3.2}) and (\ref{eqCor3.1.2}) seem more powerful than (\ref{eqth1.1}) and (\ref{eqth2.1}) because the process for approximation is independent of other random variables considered. (\ref{eqth1.1}) and (\ref{eqth2.1}) are helpful for establishing the strong convergence such as the law of the iterated logarithm.

It is of interest to find the  distribution of $\psi_{\infty}$. In the case that  $Y_{0,k}$ and  $U_{m,k}$, $k=1,2$, are all integers,  in a recent manuscript of Zhang, et al (2010) it is proved that the distribution of $\psi_{\infty}$ is absolutely continuous  and is determined by $Y_{0,1}$, $Y_{0,2}$  and the distributions of $U_{1,1}$ and $U_{1,2}$, if
$\ep[U_{1,k}\log^+U_{1,k}]<\infty$, $k=1,2$.  In the general case, the distribution of $\psi_{\infty}$ is still a open problem. In our Corollaries \ref{cor3}  and \ref{cor3.1}, by applying the Gaussian approximation and a clever idea of Aletti, May and Secchi (2009) we show that $\psi_{\infty}$ has no point masses in $[0,\infty]$ under the assumption of finite $(2+\epsilon)$-th moments. The next step is to show that the distribution is absolutely continuous. Unfortunately, as discussed in Aletti, May and Secchi (2009), the almost-sure conditional central limit theorems (\ref{eqCor2.3.1}) and (\ref{eqCor3.1.3}) are not enough to prove the absolute continuity. Our method and that of Aletti, May and Secchi (2009) depends on the martingale approach, which is not a very powerful tool to derive the limit distribution which is not normal. To find the exact distribution of $\psi_{\infty}$ needs new methods.

For a special case that  $\pr(U_{m,k}=0)=p$ and $\pr(U_{m,k}=\alpha)=1-p$ with $\alpha>0$ and $0<p\le 1$, in Aletti, May and Secchi (2007) it is shown that  the distribution of $Z_{\infty}$ is a beta distribution, and hence the probability that $\psi_{\infty}=Z_{\infty}/(1-Z_{\infty})$ falls into any subset of $[0,\infty]$ will not be zero if the Lebesgue measure of this subset is positive. So,  it is also of interest to prove in the general case that the probability of $\psi_{\infty}$ falling into any nonempty subinterval  of $(0,\infty)$ will not be zero.

Finally, this paper only consider the two-color urn model. In the manuscript of Zhang, et al (2010), the asymptotic properties for a multi-color reinforced urn model are studied. It is expected to approximate the urn components after being suitably normalized by a multi-dimensional Gaussian process.  The Skorokhod embedding method used in this paper does not work for the multi-dimension case. Though strong approximations for multi-dimensional martingales can be found in literature, for example,  Monrad and Philipp (1991), Eberlein (1986) and Zhang (2004),  the  martingales concerning to the reinforced urn model usually do not satisfied a condition  that  the asymptotic conditional variability  is   $\mathcal{F}_k$-measurable for some fixed $k$   (c.f., (\ref{eqVarApp})), which is needed in  the approximation theorems for multi-dimensional martingales. A new approach is needed for approximating the multi-color reinforced urn models.

 %%%%%%%%%%%%%%%%%%%%%%%%%%%%%%%%%%%%%%%%%%%%%%%%%%%%%%%%%%%%%%%%%%%
\appendix
\section{Proof of the main results}\label{sectionproofs}
\setcounter{equation}{0}
Recall $|\bm Y_n|=Y_{n,1}+Y_{n,2}$, $\pr(X_{n,k}=1|\mathcal{F}_{n-1})=Y_{n-1,k}/|\bm Y_{n-1}|$.
We first prove Theorems \ref{th01} and \ref{th02}.

 {\it Proofs of Theorems \ref{th01} and \ref{th02}.} Without loss of generality, assume $1<q\le 2$. It is obvious that Theorem \ref{th01} follows from Theorem \ref{th02} with  $Z_{\infty}=\psi_{\infty}/(1+\psi_{\infty})$. For Theorem \ref{th02},  let
  $$ Q_n = \frac{1}{m_1}\log Y_{n,1}-\frac{1}{m_2}\log Y_{n,2}. $$
  Then $Q_n=\frac{1}{m_1}\log \psi_n$. So, it is sufficient to show that $Q_n$ converges almost sure to a finite limit. Write
  \begin{align} \label{eqDeltaQ}\Delta Q_n =&Q_n-Q_{n-1}\nonumber\\
  &=\frac{1}{m_1} X_{n,1} \log\left(1+\frac{U_{n,1}}{Y_{n-1,1}}\right)-
  \frac{1}{m_2} X_{n,2} \log\left(1+\frac{U_{n,2}}{Y_{n-1,2}}\right)\nonumber\\
  =&\left[ X_{n,1}\frac{U_{n,1}/m_1}{Y_{n-1,1}}-X_{n,2}\frac{U_{n,2}/m_2}{Y_{n-1,2}}\right]\nonumber\\
  &+\left[-\frac{1}{m_1} X_{n,1} f\Big(\frac{U_{n,1}}{Y_{n-1,1}}\Big)+
  \frac{1}{m_2} X_{n,2} f\Big(\frac{U_{n,2}}{Y_{n-1,2}}\Big)\right]\\
  :=& \Delta Q_n^{(1)}+\Delta Q_n^{(2)}\nonumber  ,
\end{align}
where
$f(x)=x-\log(1+x)$ satisfying $0\le f(x)\le x^q$ for $x\ge 0$. We need to prove the almost sure convergence of the random series $\sum_{n=1}^{\infty}\Delta Q_n$. We first consider the terms in the second bracket above. Denote $\mathscr{A}_n=\sigma(\mathcal{F}_n, X_{n+1,1},X_{n+1,2})$. Note that $Y_{n,k}\approx N_{n,k}$ a.s. by Lemma A.4 of Hu and Zhang (2004), and then $|\bm Y_n|\approx N_{n,1}+N_{n,2}=n$ a.s.. We have that
\begin{align*}
&\sum_{n=1}^{\infty}\ep \left[X_{n,k} f\Big(\frac{U_{n,k}}{Y_{n-1,k}}\Big)\Big|\mathscr{A}_{n-1}\right]
\le \sum_{n=1}^{\infty}X_{n,k} \ep\left[\Big(\frac{U_{n,k}}{Y_{n-1,k}}\Big)^q\Big|\mathscr{A}_{n-1}\right]\\
&\le \sum_{n=1}^{\infty}X_{n,k}  \frac{\ep U_{n,k}^q}{Y_{n-1,k}^q}\le C\sum_{n=1}^{\infty} \frac{X_{n,k}}{(1+N_{n,k})^q}\\
&\le C\sum_{n=1}^{\infty} \int_{N_{n-1,k}}^{N_{n,k}}\frac{1}{(1+x)^q} dx \le C\int_0^{\infty}\frac{dx}{(1+x)^q}<\infty \;\; a.s.
\end{align*}
Hence $\sum_{n=1}^{\infty}\Delta Q_n^{(2)}$ converges almost surely. For $\{\Delta Q_n^{(1)}\}$, it is easily seen that it is a sequence of martingale differences with respect to the $\sigma$-filtration $\{\mathcal{F}_n\}$, and
$$ \ep\big[\big|\Delta Q_n^{(1)}\big|^q|\mathcal{F}_{n-1}\big]\le \frac{\ep[|U_{1,1}/m_1|^q]}{|\bm Y_{n-1}|Y_{n-1,1}^{q-1}}+\frac{\ep[|U_{1,2}/m_1|^q]}{|\bm Y_{n-1}|Y_{n-1,2}^{q-1}}. $$
So, it is sufficient to show that $\sum_{n=1}^{\infty}1/(|\bm Y_{n-1}|Y_{n-1,k}^{q-1})<\infty$ a.s., $k=1,2$. The proof will be completed if we have proven that $Y_{n,k}\ge   n^{\epsilon}$ for some positive $\epsilon$. Now, it is obvious that
$$ \sum_{n=1}^{\infty}(\log n)^{-q}\ep\big[\big|\Delta Q_n^{(1)}\big|^q|\mathcal{F}_{n-1}\big]\le C\sum_{n=1}^{\infty}\frac{1}{|\bm Y_{n-1}|(\log n)^q}<\infty \; a.s.,$$
which implies that
$$ \frac{1}{\log n} \sum_{l=1}^n \Delta Q_n^{(1)}\to 0\;\; a.s., \;\; \text{ and hence } \frac{Q_n}{\log n} \to 0\;\; a.s. $$
On the event $\{Y_{n,1}\ge Y_{n,2}\}$, we have that
\begin{align*}
 \log Y_{n,2}=& -m_2Q_n+\frac{m_2}{m_1}\log |\bm Y_n|+\frac{m_2}{m_1}\log \frac{Y_{n,1}}{|\bm Y_n|}\\
 \ge &
-m_2Q_n+\frac{m_2}{m_1}\log |\bm Y_n|+\frac{m_2}{m_1}\log \frac{1}{2}.
\end{align*}
On the event $\{Y_{n,1}\le Y_{n,2}\}$, we have $Y_{n,2}\ge |\bm Y_n|/2$. Note that $|\bm Y_n|\approx n$. It follows that
$$\liminf_{n\to \infty} \frac{\log Y_{n,2}}{\log n}\ge \frac{m_2}{m_1}\wedge 1\;\; a.s. $$
Similarly, we have that
$$\liminf_{n\to \infty} \frac{\log Y_{n,1}}{\log n}\ge \frac{m_1}{m_2}\wedge 1\;\; a.s. $$
The proof of Theorems  \ref{th01} and \ref{th02} is now completed. $\Box$.

\bigskip
{\it Proof of Theorem \ref{th1}.}
Let $f(x)=e^{m_1x}/(1+e^{m_1x})$. Then $f(Q_n)=Z_n$ and $f^{\prime}(Q_{\infty})=m_1Z_{\infty}(1-Z_{\infty})$.
   According to the Taylor expansion,  it is sufficient to show that $B(t)$ and $T_n$ can be constructed such that
  \begin{equation}\label{eqappforQ} Q_{\infty}-Q_n =\frac{H}{m_1}\int_{nH}^{\infty} \frac{d B(x)}{x} +o(\lambda_n)
 \;\; a.s..
\end{equation}

Recall $m_1=m_2$. It is easily shown that $|\bm Y_n|/n\to m_1$ a.s.. According to Theorem \ref{th01}, $Z_{\infty}\in (0,1)$ a.s.,
 which implies that $Y_{n,k}\approx n$ a.s., $k=1,2$.    So, for $\Delta Q_n^{(2)}$ in (\ref{eqDeltaQ}) we have
$$\sum_{l=1}^{\infty} \lambda_l^{-1}\ep[|\Delta Q_l^{(2)}|\big|\mathcal{F}_{l-1}]\le
\sum_{l=1}^{\infty} \frac{\lambda_l^{-1}}{|\bm Y_{n-1}|}\left(\frac{\sigma_1^2}{Y_{l-1,1}}+\frac{\sigma_2^2}{Y_{l-1,2}}\right)
\le C
\sum_{l=1}^{\infty} \lambda_l^{-1}l^{-2}<\infty,
$$
which implies that $\sum_{l=n+1}^{\infty}  |\Delta Q_l^{(2)}|=o(\lambda_n)$ a.s..

For $\Delta Q_n^{(1)}$, we use the truncation method. Let $\widetilde{U}_{n,k}=U_{n,k}/m_k$ $\overline{U}_{n,k}=\widetilde{U}_{n,k}I\{\widetilde{U}_{n,k}\le n^{1/p}\}$, $\overline{\sigma}_{n,k}^2=\ep \overline{U}_{n,k}^2$, $\overline{m}_{n,k}=\ep[\overline{U}_{n,k}]$,  $k=1,2$, $\overline{m}_n=\overline{m}_{n,1}-\overline{m}_{n,2}$, and
$$\Delta M_n^{(1)}=m_1n\left(X_{n,1}\frac{\overline{U}_{n,1}}{Y_{n-1,1}}-X_{n,2}\frac{\overline{U}_{n,2}}{Y_{n-1,2}}\right),$$
$$\Delta M_n=\Delta M_n^{(1)}-\ep[\Delta M_n^{(1)}|\mathcal{F}_{n-1}]=\Delta M_n^{(1)}-\frac{m_1n\overline{m}_n}{|\bm Y_{n-1}|}.$$
Then $\{\Delta M_n,\mathcal{F}_n\}$ is a sequence of martingale differences.  Note that
$$ \sum_{n=1}^{\infty}\pr(U_{n,k}/m_k> n^{1/p})\le \ep(U_{1,k}/m_k)^p<\infty, \;\; k=1,2.$$
From the Borel-Cantelli lemma, it follows that
$$ \pr(\Delta Q_n^{(1)}\ne \frac{1}{m_1n}\Delta M_n\;\; i.o.)=0. $$
Also,
\begin{align*}
& \sum_{l=n+1}^{\infty}\frac{1}{m_1l}|\ep[\Delta M_l^{(1)}|\mathcal{F}_{l-1}]|=\sum_{l=n+1}^{\infty}\frac{|\overline{m}_l|}{|\bm Y_{l-1}|}\\
\le &C\sum_{l=n+1}^{\infty}\frac{1}{l} \sum_{k=1}^2 \ep[|\widetilde{U}_{1,k}|I\{\widetilde{U}_{1,k}>l^{1/p}\}]
\le C n^{1/p-1}\sum_{k=1}^2 \ep[|\widetilde{U}_{1,k}|^p] =o(\lambda_n).
\end{align*}
Hence, we conclude that
\begin{align*}
\sum_{l=n+1}^{\infty} &\Delta Q_l =\sum_{l=n+1}^{\infty} \frac{1}{m_1 l} \Delta M_l+o(\lambda_n) \\ =&\frac{1}{m_1}\left(\sum_{l=n}^{\infty} \frac{1}{l(l+1)}M_l-\frac{M_n}{n}\right)+o(\lambda_n)\;\; a.s.
\end{align*}
For the martingale $M_n=\sum_{l=1}^n \Delta M_l$, we have
\begin{equation}\label{eqconditionE.1} \ep[(\Delta M_n)^2|\mathcal{F}_{n-1}]=\Big(\frac{m_1n}{|\bm Y_{n-1}|}\Big)^2
\Big(\frac{\ep \overline{U}_{n,1}^2}{Z_{n-1} }+\frac{\ep \overline{U}_{n,2}^2}{1-Z_{n-1} }-\overline{m}_n^2\Big),
\end{equation}
\begin{align}\label{eqconditionE.2}
\ep[|\Delta M_n|^4|\mathcal{F}_{n-1}]\le & \sum_{k=1}^2\left(\frac{m_1n}{Y_{n-1,k}}\right)^4\ep \overline{U}_{n,k}^4
 \le C(\omega)\sum_{k=1}^2 \ep \overline{U}_{n,k}^4.
\end{align}
By the Skorokhod embedding theorem (c.f., Theorem A.1 of Hall and Heyde (1980,page 269)), (possibly in an enlarged probability space) there is a standard motion $B(x)$ with a filtration $\{\mathscr{G}_n\}$ and a sequence of nonnegative stopping times  $\tau_1,\tau_2,\cdots $ with the following properties
\begin{description}
 \item[\rm (i)] $M_n=B(T_n)$, where $T_n=\sum_{i=1}^n \tau_i$;
  \item[\rm (ii)]  $\mathcal{F}_n\subset \mathscr{G}_n$, $\tau_n$ is $\mathscr{G}_n$ measurable, $\ep[\tau_n|\mathscr{G}_{n-1}]=\ep[(\Delta M_n)^2|\mathcal{F}_{n-1}]$, $\ep[\tau_n^r|\mathscr{G}_{n-1}]\le C_r\ep[(\Delta M_n)^{2r}|\mathcal{F}_{n-1}]$ for any $r\ge 1$;
  \item[\rm (iii)]  Conditional on $\mathscr{G}_n$,
  $B(T_n+x)-B(T_n)$, $x\ge 0$, is also a standard Brownian motion.
\end{description}

 Now, we verify that the Brownian motion $B(x)$ and the stopping time $T_n$ are  desirable for Property (b) and (\ref{eqappforQ}).   At first, we assume the following approximation for the conditional variance.
\begin{equation}\label{eqVarApp}
\ep[(\Delta M_n))^2|\mathcal{F}_{n-1}]=H(\omega)+o(n^{2/p-1})\;\; a.s.
\end{equation}
From (\ref{eqVarApp}) it follows that
$$ \sum_{i=1}^n \ep[\tau_n|\mathscr{G}_{i-1}]=n H(\omega)+ o(n^{2/p})\; a.s. $$
On the other hand, by (ii) and (\ref{eqconditionE.2}) we have that
$$ \sum_{n=1}^{\infty}\ep\left[\Big(\frac{\tau_n}{n^{2/p}}\Big)^2\big|\mathscr{G}_{n-1}\right]
\le C \sum_{n=1}^{\infty}\frac{\sum_{k=1}^2\ep \overline{U}_{n,k}^4}{n^{4/p}}\le C\sum_{k=1}^2\ep[(U_{1,k}/m_k)^p]<\infty.$$
By the law of large numbers of martingale, it follows that
$$ \sum_{i=1}^n (\tau_i-\ep[\tau_i|\mathscr{G}_{i-1}])=o(n^{2/p})\;\; a.s.. $$
Hence
$$ T_n=\sum_{i=1}^n \tau_i=n H+o(n^{2/p})\; a.s.. $$
Property (b) is verified.
Then it follows from the path properties of a Browian motion (c.f., Theorem 1.2.1 of Cs\"org\H o and R\'ev\'esz (1981)) that
\begin{align*}
 B(T_n)-B(nH)=o\left(\sqrt{n^{2/p}\Big(\log\frac{n}{n^{2/p}}+\log\log n}\Big)\right)=o(n\lambda_n) \;\;a.s.
\end{align*}
So
\begin{align*}
\sum_{l=n+1}^{\infty} &\Delta Q_l  =\frac{1}{m_1}\left(\sum_{l=n}^{\infty} \frac{1}{l(l+1)}B(T_l)-\frac{B(T_n)}{n}\right)+o(\lambda_n)\\
=&\frac{1}{m_1}\left(\sum_{l=n}^{\infty} \frac{1}{l(l+1)}B(lH)-\frac{B(nH)}{n}\right)+ \sum_{l=n}^{\infty} \frac{1}{l(l+1)}o(l\lambda_l)+o(\lambda_n)\\
=&\frac{1}{m_1}\left(\int_n^{\infty} \frac{B(xH)}{x^2}dx-\frac{B(nH)}{n}\right)+o(\lambda_n)
\\
=&\frac{1}{m_1}\left(H\int_{nH}^{\infty} \frac{B(x)}{x^2}dx-\frac{B(nH)}{n}\right)+o(\lambda_n)\\
=&\frac{H}{m_1}\int_{nH}^{\infty} \frac{d B(x)}{x} +o(\lambda_n) \;\; a.s..
\end{align*}

Finally, we verify (\ref{eqVarApp}). Note that $\ep [\overline{U}_{n,k}^2]\to \sigma_k^2$, $|\bm Y_n|/n\to m_1$ a.s., $k=1,2$ and  $Z_n\to Z_{\infty}$ a.s.. (\ref{eqVarApp}) is obvious for $p=2$ by (\ref{eqconditionE.1}).

For $2<p<4$, we still have   Corollary (\ref{cor1}) due to the approximation for the case of $p=2$. Hence
$$ Z_n-Z_{\infty}=O(\sqrt{n^{-1}\log\log n})=o(n^{2/p-1})\;\; a.s. $$
On the other hand,
\begin{align*}
 \frac{Y_{n,1}+Y_{n,2}}{m_1 n} -1&=\frac{\sum_{k=1}^2 \sum_{i=1}^n X_{i,k}(U_{i,k}-\ep[U_{i,k}]) }{m_1n}\\
 =&O(\sqrt{n^{-1} \log\log n })=o(n^{2/p-1}))\; a.s.
 \end{align*}
and
$$|\overline{m}_n|\le \sum_{k=1}^2\ep[U_{1,k}/m_k I\{U_{1,k}/m_k>n^{1/p}\}]=o(n^{1/p-1}), $$
$$ \sigma_k^2-\ep[\overline{U}_{n,k}^2]=\ep[(U_{1,k}/m_k)^2I\{U_{1,k}/m_k>n^{1/p}\}]=o(n^{2/p-1}), \;k=1,2. $$
(\ref{eqVarApp}) follows by (\ref{eqconditionE.1}). The proof is now completed. $\Box$

\bigskip

{\it Proof of Theorem \ref{th2}.} As before, denote $|\bm Y_n|=Y_{n,1}+Y_{n,2}$, $Q_n = \frac{1}{m_1}\log Y_{n,1}-\frac{1}{m_2}\log Y_{n,2}. $
   According to the Taylor expansion,  it is sufficient to show that $B(t)$ and $T_n$ can be constructed such that
  \begin{equation}\label{eqappforQth2} Q_{\infty}-Q_n =\frac{\sigma_1}{\sqrt{m_1m_2^{\rho}}\psi_{\infty}}
  \int_{n^{\rho}/\psi_{\infty}}^{\infty}  \frac{ d B(y)}{y^2}+ o( n^{-\rho(1+\delta)/2})
 \;\; a.s..
\end{equation}
Recall (\ref{eqDeltaQ}) and note that $Y_{n-1,2}\sim m_2 n$, $Y_{n-1,1}\sim \psi_{\infty}Y_{n,2}^{\rho}\sim \psi_{\infty}(m_2n)^{\rho}$ a.s.. It can be show that for $0<\delta_1\le 1/2$ and $\delta_1<(1-\rho)/\rho$,
$$\sum_{l=1}^{\infty}  l^{\rho(1+\delta_1)/2}\ep\Big[X_{l,k}f\Big(\frac{U_{l,k}}{Y_{l-1,k}}\Big)\Big|\mathcal{F}_{l-1}\Big]\le
\sum_{l=1}^{\infty} \frac{ l^{\rho(1+\delta_1)/2}}{l\cdot l^{\rho}}< \infty,
$$
which implies $\sum_{l=n+1}^{\infty} X_{l,k}f\big(U_{l,k}/Y_{l-1,k}\big)=o(n^{-\rho(1+\delta_1)/2})$ a.s., $k=1,2$.
Also, for the martingale differences $X_{l,2}\frac{U_{l,2}/m_2}{Y_{l-1,2}}-\frac{1}{|\bm Y_{l-1}|}$, we have
$$\sum_{l=1}^{\infty}  (l^{\rho(1+\delta_1)/2})^2\ep\Big[\Big(X_{l,2}\frac{U_{l,2}/m_2}{Y_{l-1,2}}-\frac{1}{|\bm Y_{l-1}|}\Big)^2\Big|\mathcal{F}_{l-1}\Big]\le
\sum_{l=1}^{\infty} \frac{ l^{\rho(1+\delta_1)}}{l^2}< \infty\;\;a.s.,
$$
which implies $\sum_{l=n+1}^{\infty}\big(X_{l,2}\frac{U_{l,2}/m_2}{Y_{l-1,2}}-\frac{1}{|\bm Y_{l-1}|}\big)=o(n^{-\rho(1+\delta_1)/2})$ a.s.. Similarly, we can show that
 $\sum_{l=n+1}^{\infty}\big(X_{l,1}\frac{U_{l,1}/m_2}{Y_{l-1,1}}-\frac{1}{|\bm Y_{l-1}|}\big)=o(n^{-\rho /2}\log n)$ a.s.
It follows that
 \begin{equation}\label{eqproofth2.5}
 Q_{\infty}-Q_n =\sum_{l=n+1}^{\infty}\left(X_{l,1}\frac{U_{l,1}/m_1}{Y_{l-1,1}}-\frac{1}{|\bm Y_{l-1}|}\right)+o(n^{-\rho(1+\delta_1)/2})\; a.s.
 \end{equation}
 and
  \begin{equation}\label{eqproofth2.6}
 Q_{\infty}-Q_n = o(n^{-\rho/2}\log n)\; a.s.
 \end{equation}
Define
\begin{equation}\label{eqproofth2.7}  \Delta M_n =\frac{\sqrt{\rho m_2}(m_2n)^{\rho/2}n^{\rho/2}}{\sigma_1}\left(X_{n,1}\frac{U_{n,1}/m_1}{Y_{n-1,1}}-\frac{1}{|\bm Y_{n-1}|}\right).
\end{equation}
Then
 \begin{align*}
 \ep[(\Delta M_n)^2|\mathcal{F}_{n-1}]=&\frac{\rho m_2(m_2n)^{\rho}n^{\rho}}{|\bm Y_{n-1}|Y_{n-1,1}}-\frac{\rho m_2(m_2n)^{\rho}n^{\rho}}{\sigma_1^2|\bm Y_{n-1}|^2}.
 \end{align*}
 Next, we first show that
  \begin{equation}\label{eqproofth2.8}
 \ep[(\Delta M_n)^2|\mathcal{F}_{n-1}]
 = \frac{\rho}{\psi_{\infty}}n^{\rho-1}(1+o(n^{-\rho\delta_1}) )\;\; a.s.
 \end{equation}
From (\ref{eqproofth2.6}) and the Taylor expansion, we conclude that
$$\frac{Y_{n,1}}{(Y_{n,2})^{\rho}}-\psi_{\infty}=o(n^{-\rho/2}\log n)  \;\; a.s. $$
 On the other hand,
 \begin{align*}
  \frac{ Y_{n,2}}{m_2n}=&1-\frac{1}{m_1}\frac{Y_{n,1}}{n}
 +\sum_{k=1}^2\frac{\sum_{l=1}^n X_{l,k}(U_{l,k}/m_k-\ep[U_{l,k}/m_k])}{n}\\
 =& 1-\frac{1}{m_1}\frac{\psi_{\infty}(Y_{n,2})^{\rho}}{n}+o(n^{-\rho/2-1}\log n)+O(n^{-1/2}(\log\log n)^{1/2})\\
 =&1-O(n^{\rho-1})+o(n^{-\rho\delta_1})=1+o(n^{-\rho\delta_1})\;\; a.s.
 \end{align*}
It follows that
$$\frac{|\bm Y_n|}{m_2n} =\frac{ Y_{n,2}}{m_2n}+\frac{ Y_{n,1}}{m_2n}=1+o(n^{-\rho\delta_1})\; \text{and }\; \frac{Y_{n,1}}{(m_2n)^{\rho}}-\psi_{\infty}=o(n^{-\rho\delta_1}) \; a.s. $$
 (\ref{eqproofth2.8}) is verified. From (\ref{eqproofth2.8}), it follows that
 $$\sum_{l=1}^n \ep[(\Delta M_n)^2|\mathcal{F}_{n-1}]=\frac{n^{\rho}}{\psi_{\infty}}(1+o(n^{-\rho\delta_1}) )\;\; a.s. $$
 On the other hand, for $0<\delta_2<1/2-1/p$,
 \begin{align*}
  \sum_{n=1}^{\infty}\frac{\ep[|\Delta M_n|^p|\mathcal{F}_{n-1}]}{(n^{\rho(1-\delta_2)})^{p/2}}
 \le& C \sum_{n=1}^{\infty} \frac{n^{\rho p}}{n^{p(1-\delta_2)\rho/2}}\frac{Y_{n-1,1}}{Y_{n-1,1}^p|\bm Y_{n-1}|}\\
 \le &C\sum_{n=1}^{\infty}\frac{1}{n^{p(1-\delta_2)\rho/2}}\frac{1}{n^{1-\rho}}<\infty.
 \end{align*}
 So, similarly as in the proof of Theorem \ref{th1}, by the Skorokhod embedding theorem,  the  standard motion $B(x)$, the filtration $\{\mathscr{G}_n\}$ and the stopping times $\{T_n\}$ can be   constructed such that
 $M_n=B(T_n)$ and
 $$ T_n=\sum_{l=1}^n \ep[(\Delta M_l)^2|\mathcal{F}_{l-1}]+o(n^{\rho(1-\delta_2)})
 =\frac{1}{\psi_{\infty}}n^{\rho} +o(n^{\rho(1-\delta_1\wedge \delta_2)} ). $$
Denote $\delta_0=\min\{(1-\rho)/\rho,1/2-1/p\}$. It is remained to verify (\ref{eqth2.1}). By the Properties (b$^{\prime}$) and the path properties of a Brownian motion, we have  for any $0<\delta<\delta_0$,
 $$ M_n-B(n^{\rho}/\psi_{\infty})=o(n^{\rho(1-\delta)/2}) \;a.s. $$
 Hence
 \begin{align*}
 &\sum_{l=n+1}^{\infty}\frac{\Delta M_l}{l^{\rho}} =\sum_{l=n}^{\infty}\left(\frac{1}{l^{\rho}}-\frac{1}{(l+1)^{\rho}}\right)M_l-\frac{M_n}{n^{\rho}}\\
 =&\sum_{l=n}^{\infty}\left(\frac{1}{l^{\rho}}-\frac{1}{(l+1)^{\rho}}\right)B(l^{\rho}/\psi_{\infty})
 -\frac{B(n^{\rho}/\psi_{\infty})}{n^{\rho}}\\
 \qquad & +\sum_{l=n}^{\infty}\frac{o(l^{\rho(1-\delta)/2})}{l^{1+\rho}}+\frac{o(n^{\rho(1-\delta)/2})}{n^{\rho}}\\
 =& \int_n^{\infty} \frac{\rho B(x^{\rho}/\psi)}{x^{1+\rho}}dx -\frac{B(n^{\rho}/\psi_{\infty})}{n^{\rho}}+
 o(n^{-\rho(1+\delta)/2})
\\
 =& \frac{1}{\psi_{\infty}} \int_{n^{\rho}/\psi_{\infty}}^{\infty} \frac{B(x )}{x^2}dx -\frac{B(n^{\rho}/\psi_{\infty})}{n^{\rho}}+
 o(n^{-\rho(1+\delta)/2})
\\
 =& \frac{1}{\psi_{\infty}} \int_{n^{\rho}/\psi_{\infty}}^{\infty} \frac{d B(x )}{x}dx +
 o(n^{-\rho(1+\delta)/2})  \;\; a.s.
 \end{align*}
(\ref{eqappforQth2}) is now proved by noting that (\ref{eqproofth2.5}), (\ref{eqproofth2.6}) and $\rho m_2=m_1$.
And hence (\ref{eqth2.1}) is verified. $\Box$

\begin{rem}
Using the truncation method as in the proof of Theorem \ref{th1}, we can proved that (\ref{eqth2.1}) remains true under the assumption of only finite second moments if $n^{-\rho(1+\delta)/2}$ is replaced by $n^{-\rho/2}(\log\log n)^{1/2}$. This implies that the law of iterated logarithm  (\ref{eqCor3.1.1}) remains true when $\ep U_{1,k}^2<\infty$, $k=1,2$.
\end{rem}
%\end{document}


\bigskip


\begin{thebibliography}{99}
\footnotesize



\bibitem{AMS07}
{\sc Aletti, G., May, C.} and {\sc Secchi, P.} (2007). On the distribution of the limit proportion for a two-color, randomly reinforced urn with equal reinforcement distributions. {\em Adv. Appl. Probab.}, {\bf 39}: 690-707.


\bibitem{AMS09}
{\sc Aletti, G., May, C.} and {\sc Secchi, P.} (2009). A central limit theorem, and related results, for a two-color randomly reinforced urn. {\em  Adv. Appl. Probab.},  {\bf 41}: 829-844.



\bibitem{BH05}
{\sc Bai, Z. D.} and {\sc Hu, F.} (2005). Strong consistency and
asymptotic normality for urn models. {\em Ann. Appl. Probab.}, {\bf
12}: 914-940.

\bibitem{BHR02}
 {\sc Bai, Z. D.}, {\sc Hu, F.}  and {\sc Rosenberger, W. F.}  (2002).
 Asymptotic properties of adaptive
designs for clinical trials with delayed response. {\em Ann.
Statist.},  {\bf 30}: 122-139.

\bibitem{BHZ02}
{\sc Bai, Z. D.}, {\sc Hu, F.} and {\sc Zhang, L. X.}  (2002). The Gaussian approximation theorems for urn models and their applications. {\em Ann. Appl.Probab.}, {\bf 12}: 1149-1173.

\bibitem{Beggs05}
{\sc Beggs, A. W.} (2005). On the convergence of reinforcement
learning. {\em J.  Econom. Theory} , {\bf 122}(1): 1-36.

\bibitem{CPR09}
{\sc Chauvin, B., Pouyanne, N.} and {Sahnoun, R.} (2009). Limit distributions for large P\'olya urns.
http://arxiv.org/abs/0907.1477

\bibitem{Cri09}
{\sc Crimaldi, I.} (2009).
An almost sure conditional convergence result
and an application to a generalized P\'olya Urn. {\em International Mathematical Forum}, {\bf 4} (23): 1139--1156

\bibitem{DFL98}
 {\sc Durham, S. D.}, {\sc Flournoy, N.} and {\sc Li, W.} (1998).
A sequential design for maximizing the probability of a favourable
response. {\em Canad. J. Statist.}, {\bf 26} (3):  479-495.

\bibitem{DY90}
{\sc Durham, S. D.} and {\sc Yu, K. F.}  (1990).
Randomized play-the leader rules for sequential sampling from two
populations. {\em Probability in Enginerring and Information
Science}, {\bf 26} (4): 355-367.


\bibitem{Eber86}
{\sc Eberlein, E.} (1986). On strong invariance principles under dependence. {\em Ann. Probab.}, {\bf 14}: 260每270.


\bibitem{EggP23}
{\sc Eggenberger, F. } and {\sc   P∩olya, G.} (1923). Uber die Statistik verketteter Vorg\"ange. {\em Zeitschrift 	 Angew. Math. Mech.}, {\bf 3}: 279每289.

\bibitem{ErevRoth98}
{\sc Erev, I.} and {\sc Roth, A.} (1998). Predicting how people play
games: reinforcement learning in experimental games with unique,
mixed strategy equilibria.  {\em Amer. Econ. Rev.}, {\bf 88}:
848-881.

\bibitem{HallHeyde80}
{\sc Hall, P.} and {\sc Heyde, C. C.} (1980). {\it Martingale Limit
Theory and its Applications}. Academic Press, London.

\bibitem{HopkinsPosch05}
{\sc Hopkins, E. } and {\sc Posch, M.} (2005). Attainability of boundary points under reinforcement learning. {\em Games Econom. Behavior}, {\bf 53}: 110-125.

\bibitem{HansonRusso83}
{\sc Hanson, D. L.} and {\sc  Russo, Ralph P.}
 (1983). Some results on increments of the Wiener process with applications to lag sums of i.i.d.r.v. {\em Ann. Probab.}, {\bf 11}: 609-623.

\bibitem{HR06}
{\sc Hu, F.} and {\sc Rosenberger, W. F.} (2006).
{\em The Theory of Response-Adaptive Randomization in Clinical
Trials}, John Wiley and Sons, Inc., New York.



\bibitem{HZ04}
{\sc Hu, F. \& Zhang, L.-X.}  (2004).  Asymptotic properties of
doubly adaptive biased coin designs for multi-treatment clinical
trials.  {\em Ann. Statist.}, {\bf 32}: 268-301.


\bibitem{Janson04}
{\sc Janson, S.} (2004). Functional limit theorems for multitype
branching processes and generalized P\'olya urns. {\em Stochastic
Process. Appl.}, {\bf 110}: 177-245.

\bibitem{Janson06}
{\sc Janson, S.} (2006). Limit theorems for triangular urn schemes. {\em Probab. Theory Relat. Fields}, {\bf 134}:417每452.

\bibitem{LDF96}
 {\sc Li, W.}, {\sc Durham, S. D.} and  {\sc Flournoy, N.}
(1996).
 Randomized polya urn designs. {\em Proceedings of the Biometric
 Section of the Statistical Association}: 166-170.

\bibitem{MF09}
{\sc May, C. } and {\sc Flournoy, N.} (2009). Asymptotics in response-adaptive designs generated by a two-color, randomly reinforced urn.
{\em Ann. Statist.},  {\bf 37}(2): 1058-1078.


\bibitem{MH02}
{\sc Martin, C. F.} and {\sc Ho, Y. C.} (2002). Value of information in
the Polya urn process. {\em Information Sciences}, {\bf 147}: 65-90.

\bibitem{MP00}
{\sc Melfi, V. F.} and  {\sc Page, C.}  (2000).  Estimation after adaptive
allocation.  {\em J. Statist. Plann. Inf.},  {\bf 87}: 353--363.


\bibitem{MonPhill91}
{\sc Monrad, D.} and {\sc Philipp,W.} (1991).  Nearby variables with nearby conditional laws and a strong approximation theorem
for Hilbert space valued martingales. {\em Probab. Theory Relat. Fields}, {\bf 88}: 381每-404.




\bibitem{MPS06a}
{\sc Muliere, P.}, {\sc Paganoni, A. M.} and {\sc Secchi, P.} (2006a).
A randomly reinforced urns. {\em J. Statisit. Plann. Inference}, {\bf 136}(6):1853-1874.

\bibitem{MPS06b}
{\sc Muliere, P.}, {\sc Paganoni, A. M.} and {\sc Secchi, P.} (2006b).
Randomly reinforced urns for clinical trials with
continuous responses. In {\em SIS〞Proceedings of the XLIII Scientific Meeting}, 403每
414. Cleup, Padova.



\bibitem{PS07}
{\sc Paganoni, A.} and {\sc Secchi, P.} (2007). A numerical study
for comparing two response-adaptive designs for continuous treatment
effects. {\em Statstical Methods and Applications}, {\bf 16}: 321每346.

\bibitem{Polya31}
{\sc P\'olya, G.} (1931). Sur quelques points de la th\'eorie des probabilit\'es. {\em Ann. Inst. Poincar\'e}, {\bf  1}: 117每161.

\bibitem{Z04}
{\sc Zhang, L. X.} (2004). Strong approximations of martingale vectors and its applications in Markov-Chain adaptive designs. {\sc Acta Math. Appl. Sinica, English Series}, {\bf 20}(2): 337--352.


\bibitem{ZH09}
{\sc Zhang, L. X.} and {\sc Hu, F.} (2009). The Gaussian approximation for multi-color generalized Friedman's urn model. {\em Science in China, Ser. A}, {\bf 52} (6): 1305-1326.

\bibitem{ZHC06}
{\sc Zhang, L. X.}, {\sc Hu, F.} and {\sc Cheung, S. H.} (2006). Asymptotic theorems of sequential estimation-adjusted urn models for clinical trials. {\em Ann. Appl. Probab.}, {\bf 16}(1): 340-369

\bibitem{ZHCC10}
{\sc Zhang, L. X.}, {\sc Hu, F.}, {\sc Cheung, S. H.} and {\sc Chan, W. S.} (2010).
Asymptotic properties of multi-color randomly reinforced P\'olya urns. {\em Manuscript.}
\end{thebibliography}
\end{document}